\documentclass[12pt]{article}
\usepackage{amsmath}
\usepackage{amssymb}
\usepackage{amsfonts}
\usepackage{amsthm}
\usepackage{enumerate}
\usepackage{url}

\def\init{\setcounter{equation}{0}}

\newtheorem{theorem}{Theorem}[section]

\newtheorem{corollary}[theorem]{Corollary}

\def\Res{{\operatorname{Res}}}

\def\d{{\operatorname{d}}}

\def\Re{\hbox{Re}\,}

\def\ii{{\mathrm{i}}}
\def\e{{\mathrm{e}}}

\numberwithin{equation}{section}
\newenvironment{acknowledgements}{\noindent{\bf Acknowledgements}\bigskip}{}

\begin{document}

\title{A note on a curious formula for Euler's constant}
\author{Mathew D. Rogers\\
\small{\textit{Department of Mathematics, University of British
        Columbia}}\\
        \small{\textit{Vancouver, BC, V6T-1Z2, Canada}}\\
        \small{\textit{matrogers@math.ubc.ca}}}
        \maketitle

\abstract{In this short note we will use the residue theorem to
establish a formula for Euler's constant.  In particular, we offer a
slightly generalized version of an interesting infinite series due
to Flajolet, Gourdon, and Dumas.}
\section{Properties of $\gamma$}
\label{intro} \init
Recall that $\gamma$ is usually defined in terms of the partial sums
of the harmonic series:
\begin{equation}\label{gamma definition}
\gamma:=\lim_{n\rightarrow\infty}\left(\sum_{j=1}^{n}\frac{1}{j}-\log
n\right)=0.577215664901\dots.
\end{equation}
While this definition is interesting, it is extremely inefficient
for actually calculating the constant numerically.  Since Euler's
constant arises in many mathematical contexts \cite{We}, and since
it is often considered to be the most important ``special constant"
after $\pi$ and $e$, it is naturally an interesting problem to find
rapidly converging series expansions for $\gamma$.

While $\pi$ and $e$ were both proved to be irrational in the
eighteenth century, the arithmetic nature of $\gamma$ remains a
mystery. Recall that Fourier obtained the irrationality of $e$ as a
simple consequence of Euler's famous infinite series:
\begin{equation}\label{e definition}
e=\sum_{k=0}^{\infty}\frac{1}{k!}.
\end{equation}
While it seems virtually certain that no identities like \eqref{e
definition} exist for $\gamma$, many interesting formulas are
scattered throughout the literature.  For example, Sondow and
Zudilin proved several rational expansions for $\gamma$ in
\cite{SZ}. A variety of series expansions involving the Riemann zeta
function, and multiple identities containing logarithms also exist
(see \cite{Go} or \cite{We}). In the following theorem we offer a
slightly generalized version of a formula due to Flajolet, Gourdon
and Dumas (see page $28$ in \cite{Fl}).

\begin{theorem} Suppose that $x>0$ and $w> 0$, then the
following identity is true:
\begin{equation}\label{gamma identity}
\begin{split}
\gamma=&\frac{w}{2}-\log(x)-w\sum_{k=0}^{\infty}e^{-x e^{w
k}}+w\sum_{k=1}^{\infty}\frac{(-1)^{k+1}}{k!}\frac{x^k}{e^{w k}-1}\\
&+\sum_{\substack{k=-\infty\\k\not=0}}^{\infty}\Gamma\left(\frac{2\pi\ii
k}{w}\right)x^{-2\pi\ii k/w}.
\end{split}
\end{equation}
\end{theorem}
\begin{proof} Let $\phi_w(x)$ be defined by
\begin{equation*}
\phi_w(x):=\sum_{k=0}^{\infty}e^{-x e^{w k}},
\end{equation*}
and suppose that $w>0$.  Since this sum converges uniformly it is
easy to calculate the Mellin transform of $\phi_w(x)$:
\begin{equation*}
\int_{0}^{\infty}x^{a-1}\phi_w(x)\d x=\frac{\Gamma(a)}{1-e^{-a w}}.
\end{equation*}
Inverting the Mellin integral yields
\begin{equation}\label{phi mellin transform}
\phi_w(x)=\int_{C}\frac{\Gamma(a)x^{-a}}{1-e^{-a w}}\d a,
\end{equation}
where ${C}$ is a closed contour which runs vertically along the line
$(\frac{1}{2}-\ii\infty,\frac{1}{2}+\ii\infty)$, and then encircles
the negative half plane. Notice that the integrand in \eqref{phi
mellin transform} has simple poles on the imaginary axis at $a=2\pi
\ii k/w$ for $k\in\{\pm 1,\pm 2,\dots\}$, simple poles at the
negative integers, and a double pole at $a=0$. It follows easily
from the residue theorem that
\begin{align*}
\phi_w(x)=&\mathop{\Res}_{a=0}\left(\frac{\Gamma(a)x^{-a}}{1-e^{-a
w}}\right)+\sum_{\substack{k=-\infty\\
k\not=0}}^{\infty}\mathop{\Res}_{a=2\pi\ii
k/w}\left(\frac{\Gamma(a)x^{-a}}{1-e^{-a
w}}\right)+\sum_{k=1}^{\infty}\mathop{\Res}_{a=-k}\left(\frac{\Gamma(a)x^{-a}}{1-e^{-a
w}}\right)\notag\\
=&\frac{1}{2}-\frac{\gamma}{w}-\frac{\log(x)}{w}+\frac{1}{w}\sum_{\substack{k=-\infty\\k\not=0}}^{\infty}\Gamma\left(\frac{2\pi\ii
k}{w}\right)x^{-2\pi\ii
k/w}+\sum_{k=1}^{\infty}\frac{(-1)^{k+1}}{k!}\frac{x^k}{e^{w k}-1}.
\end{align*}
Rearranging the last identity completes the proof. $\blacksquare$
\end{proof}

Notice that a wide variety of formulas involving $\gamma$, $\pi$,
and values of the Riemann zeta function at odd integers can be
obtained by modifying the definition of $\phi_w(x)$ (for instance by
adding arithmetic functions to the sum).  As a simple example we can
show that
\begin{equation*}
6\gamma^2+\pi^2=1+12\sum_{k=1}^{\infty}k\left(e^{-e^{k}}-e^{-e^{-k}}+1\right)\\
-24\sum_{k=1}^{\infty}\Re\left[\Gamma'(2\pi \ii k)\right].
\end{equation*}
 We will also point out that equation \eqref{gamma
identity} reduces to a classical formula involving the exponential
integral when $w$ approaches zero \cite{Go}:
\begin{equation*}
\begin{split}
\gamma=&-\log(x)-\int_{0}^{\infty}e^{-x e^{u}}\d
u+\sum_{k=1}^{\infty}\frac{(-1)^{k+1}}{k}\frac{x^k}{k!}\\
 =&-\log(x)-\int_{x}^{\infty}\frac{e^{-u}}{u}\d
 u+\int_{0}^{x}\frac{1-e^{-u}}{u}\d u
\end{split}
\end{equation*}
To prove this limiting case, simply observe that the vanishing of
the Gamma sum follows immediately from a standard estimate
\cite{Gr}:
\begin{equation*}
|\Gamma\left(2\pi \ii k/w\right)|=O\left(e^{-\pi^2 k/w}\right),
\end{equation*}
as $k\rightarrow\infty$. Finally, we can recover the formula of
Flajolet, Gourdon and Dumas \cite{Fl} by setting $w=\ln(2)$ in
equation \eqref{gamma identity}.

\begin{corollary}The following identities are true:
\begin{align}
\gamma=&-2\sum_{k=0}^{\infty}e^{-e^{2k+1}}+\sum_{k=1}^{\infty}\frac{(-1)^{k+1}}{k!}\frac{1}{\sinh(k)}-2\sum_{k=1}^{\infty}(-1)^{k+1}\Re\left[\Gamma(\pi
\ii k)\right]\label{formula 1},\\
\gamma=&\frac{1}{2}-\sum_{k=0}^{\infty}e^{-e^k}+2\sum_{k=1}^{\infty}\frac{(-1)^{k+1}}{k!}\frac{1}{e^k-1}+2\sum_{k=1}^{\infty}\mathop{Re}\left[\Gamma(2\pi\ii
k)\right]\label{formula 2}.
\end{align}
\end{corollary}
\begin{proof} Equation \eqref{formula 1} follows from setting $x=e$ and $w=2$
in \eqref{gamma identity}, while \eqref{formula 2} follows from
setting $x=w=1$.$\blacksquare$
\end{proof}

%
Despite converging rapidly (notice that $\Re\left[\Gamma(\pi \ii
k)\right]=O(\e^{-\pi^2 k/2})$ as $k\rightarrow\infty$) both
equations \eqref{formula 1} and \eqref{formula 2} suffer from a
serious computational drawback. In particular, there is no obvious
way to calculate $\Gamma(\pi \ii k)$ recursively with respect to
$k$.  This fact most likely precludes the possibility of using
either \eqref{formula 1} or \eqref{formula 2} to obtain more than a
few thousand decimal digits of $\gamma$. As of December 2006, the
value of $\gamma$ was known to more than $116$ million digits
\cite{Yee}. Although a computational record may therefore seem far
out of reach, I will point out that the first sum in equation
\eqref{formula 1} would only need around $9$ terms to reach that
level of accuracy (the tenth term is $e^{-e^{21}}\approx
7\times10^{-572,754,397}\ll 10^{-116,000,000}$). Thus, it may be a
worthwhile endeavor to try to accelerate the convergence of the the
second and third sums in both equations \eqref{formula 1} and
\eqref{formula 2}.
\bigskip

\begin{acknowledgements}

The author thanks David Boyd, Jonathan Sondow, and Wadim Zudilin for
the useful discussions and encouragement.  Additionally, the author
thanks Jonathan Sondow for offering several useful corrections and
references.
\end{acknowledgements}


\end{document}